\documentclass[11pt,twoside,a4paper,notitlepage]{article}

\usepackage{graphicx}
\usepackage{mathtools}
\usepackage{enumerate}
\usepackage{amssymb}
\usepackage{amsthm}
\usepackage{mathptmx} 
\usepackage{csquotes}
\usepackage{latexsym}
\usepackage{tikz}
\usepackage[utf8]{inputenc}
\usepackage{graphicx}
\usepackage[british]{babel}

\def\.#1{\hfill#1\kern.5em\vrule\kern-.5em}

\DeclarePairedDelimiter\floor{\lfloor}{\rfloor}
\theoremstyle{definition}
\newtheorem{theorem}{Theorem}

\newtheorem{definition}{Definition}
\newtheorem{example}{Example}
\newtheorem{remark}{Remark}

\begin{document}

\title{Sizes of Countable Sets}

\author{Kate\v{r}ina Trlifajov\'{a}}


\date{}

\maketitle

\begin{abstract}

This paper introduces the notion of size of countable sets that preserves the Part-Whole Principle and generalizes the notion of cardinality of finite sets. The sizes of natural numbers, integers, rational numbers and all their subsets, unions and Cartesian products are algorithmically enumerable up to one element as sequences of natural numbers. 
The method is similar to that of \emph{Theory of Numerosities} (Benci \& Di Nasso 2019) but in comparison, it is motivated by Bolzano's concept of infinite series from his \emph{Paradoxes of the Infinite} (Bolzano 1851/2004), it is constructive because it does not use ultrafilters, and set sizes are uniquely determined. The results mostly agree with those of Theory of Numerosities, but some differ, such as the size of rational numbers. However, set sizes are just partially and not linearly ordered. \emph{Quid pro quo.}

\end{abstract}

\section{Introduction}

The introduction of actual infinity in mathematics has naturally led to two questions. The first concerns the existence of actually infinite collections, while the second deals with the way in which the actual infinite collections could be compared, for instance (Mancosu 2009, p. 612).  Two great defenders of actual infinity, Bernard Bolzano and Georg Cantor, answered the first question positively. 
However, their answers to the second question differ. When comparing infinite sets, one of two mutually exclusive principles must be chosen (Parker 2009, p. 93):

\begin{enumerate}
\item \emph{The Part-Whole Principle} (PW): \enquote{The whole is greater than its part.}\footnote{This is the 5\textsuperscript{th} Common Notion from Euclid's Element.} 
\item \emph{Hume's Principle} (HP): \enquote{Two sets have the same size if and only if there is a one-to-one correspondence between them.}\footnote{Frege states this principle with reference to Hume's \emph{Treatise of Human Nature} Book 1, Part 3, Section 1 in \emph{Grundlagen der Arithmetik}, §63. 
HP plays a central role in Frege's philosophy of mathematics. The problematic Law V, which is the source of inconsistency, is employed in his \emph{Die Grundlagen} only in the proof of HP. It is often used in \emph{Grundgesetze}, but it seems that Frege would probably have achieved the same results using HP alone. \enquote{Modulo uses of value-ranges that are essentially just for
convenience, Part II of \emph{Grundgesetze} really does contain a formal derivation
of axioms for arithmetic from HP.} (Heck 2019, p. 3 - 4)}. 
\end{enumerate}
The inconsistency of these two principles for infinite collections is well illustrated by Galileo's famous paradox of reflexivity.


\subsection{Galileo's Paradox Revisited}

 Galileo compares the infinite collection of all natural numbers and the infinite collection of their squares; that is, the second powers of natural numbers. 
$$1, \quad 2, \quad 3, \quad \ 4, \quad \ 5, \ \dots$$
$$\ \ 1, \quad 4, \quad 9, \quad 16, \quad 25,\ \dots $$
On the one hand, Galileo explains that there are more numbers, including squares and non-squares, than the squares alone. On the other hand, he also explains that there are as many squares as their roots because every square has its own root and every root its own square, and all numbers are roots. Galileo resolves this contradiction by denying the possibility of comparing infinite collection (Galileo 1638/1914, p. 40 - 42).

The possible solution of Galileo's paradox is hidden in the following question: \enquote{How do we arrive at such a collection of squares? What is its \emph{determining ground}, its \emph{way of being formed}?} There are two options. We either \emph{select} the squares from all the numbers, which makes them smaller in size. Or we \emph{create} them from numbers. In this case, we have to use Galileo's assumption that seems to be, but need not be, self-evident: \enquote{Every number is a root of some square} (Galileo 1638/1914, p. 40). Then the collection of squares has as many elements as the collection of numbers.\footnote{This concept enables Dedekind's definition of infinite sets: a set is infinite if there is a one-to-one correspondence to its proper subset.} 

Such a question was formulated by Bolzano, who argued that a mere one-to-one correspondence between two infinite collections is not a sufficient reason for their equal size, but some other reason must be added, such as that they have the same \enquote{determining ground}, see Section \ref{PW}. 

 In Zermelo-Fraenkel set theory (ZF), there are two axioms to form the set of squares of natural numbers  

\begin{enumerate}[(i)]

\item According to the \emph{Separation Schema}\footnote{Let $\varphi$ be a formula then for any set $X$ there exists a set $Y = \{u \in X;  \varphi(u)\}$ that contains
all those $u \in X$ that have property $\varphi$. (Jech 2006, p. 7).} we consider a set   
$$A = \{n \in \mathbb N; (\exists m)(m \in \mathbb N \wedge n = m^2)\}.$$ 

\item According to the  \emph{Replacement Schema}\footnote{If a class $F$ is a function, then for every set $X$ there exists a set $F(X) = \{F(u); u \in X \}$. (Jech 2006, p. 13).} we use the function $F(n) = n^2$ and get a set 
$$B  = F(\mathbb N) = \{n^2; \ n \in \mathbb N\}.$$ 

\end{enumerate}

Though the sets $A$ and $B$ are equal in ZF there is a difference in the way they are created. Following (i) alone one might conclude that $A$ is primarily a proper subset of $\mathbb N$, since  there are certainly numbers that are not squares, so there should be less squares than all numbers. However, (ii) immediately implies a one-to-one correspondence, so that $\mathbb N$ has as many elements as $B$.\footnote{The \emph{Axiom Scheme of Replacement} was added to Zermelo's 1908 axiomatisation of set theory later, in 1922. (Kanamori, 2012, p. 61). It is stronger than the \emph{Axiom Scheme of Separation} which in fact follows from \emph{Replacement} and the \emph{Axiom of Empty Set.}}

Not only does a one-to-one correspondence exist between natural numbers and their squares, but also between numbers and any of their infinite subsets. They all have the same cardinality, namely $\aleph_0$, as well as their Cartesian products, integers, rational numbers and others. Just for these sets we determine their size so that PW is valid up to one element.

\subsection{The Structure of this Paper}

The main goal is to introduce a theory of size preserving PW of some countable sets. The method is similar to that of Benci and Di Nasso's  Numerosity Theory (NT) but it differs in some significant points that are analysed in Section \ref{num}.  

The source of inspiration for this research is Bolzano's concept of the infinite. We briefly describe some principles of his theory of infinite quantities and its interpretation in contemporary mathematics in Section \ref{Bolzano}. However, the core of the theory lies in Section \ref{size} where we introduce the notion of the sizes of canonical arranged sets. It turns out that natural numbers, integers and rational numbers and their subsets can be canonically arranged, we investigate their sizes in Section \ref{numbers}. Section \ref{num} contains a comparison with the Numerosity Theory. In the last Section \ref{concl}, we answer Parker who argues that there is no good theory preserving PW.

\section{Bolzano's Concept}\label{Bolzano}

\subsection{Bolzano and the Part-Whole Principle}\label{PW}

The largest part of Bolzano's concept of the infinite can be found in his last work \emph{Paradoxes of the Infinite} (Bolzano 1851/2004) (PU) that was written in 1848. Unlike Cantor, who chose HP for a comparison of infinite sets, Bolzano did not accept one-to-one correspondence as a sufficient criterion for the equality of \textquote{pluralities} and insisted on PW even for infinite sets (Rusnock \& \v Sebest\' ik 2019, p. 536). While Cantor's theory has been accepted by the broad mathematical community, Bolzano's concept has been mostly judged from this perspective as a step in the wrong direction.\footnote{This opinion comes from Cantor, who, while praising Bolzano for his courage to defend actual infinity, nonetheless regards PU as \emph{unfounded and erroneous}. \enquote{The author lacks both the general concept
of cardinality and the precise concept of number-of-elements for a real conceptual grasp of determinate infinite numbers.} (Cantor 1883/1976, p. 78). Yet Ferreiros says about Bolzano \enquote{But after having been close to the \emph{right} point of view, he departed from it in quite a \emph{strange} direction.} The \enquote{right point of view} is the idea that cardinality is the only meaningful way to compare abstract sets with respect to the multiplicity of their elements (Ferreiros 1999, p. 75)} Only recently, some papers have challenged the hegemony of Cantor's conception, such as \emph{Measuring the Size of Infinite Collections of Natural Numbers: Was Cantor's Set Theory Inevitable?} (Mancosu 2009) or \emph{Is mathematical history written by the victors?} (Bair et al. 2013), and others. Bolzano's PU were reinvestigated and it turned out it can be well interpreted as a consistent and  meaningful (Trlifajov\'a 2017), (Bellomo \& Massas 2021).

Bolzano was aware of the existence of a one-to-one correspondence between some infinite multitudes.\footnote{Bolzano's multitude \emph{[Menge]} is the notion \enquote{sui generis} (Simons 1998, p. 87). This is a \enquote{collection of certain things or a whole consisting of certain parts} (PU \S 3) such that \enquote{the arrangement of its parts is unimportant.} (PU \S 4) Cantor later used the same German word for his sets. Although Bolzano's multitudes and Cantor's sets share many features, 
Bolzano scholars agree that Bolzano's multitudes cannot be precisely interpreted, neither as Cantor's sets nor mereologically (Rusnock 2012, p. 155).} 
He calls this property a \emph{highly remarkable peculiarity} (PU, \S 20) and warns against the assertion that one-to-one correspondence allows for the conclusion of the equal plurality of their parts. 
\begin{quote} 
Merely from this circumstance we can - as we see - in no way conclude \emph{that these multitudes are equal to one another if they are infinite} with respect to the plurality of their parts (i.e. if we disregard all differences between them). But rather they are able to have a relationship of inequality in their plurality, so that one of them can be presented as a whole, of which the other is a part. (PU \S 21).
\end{quote}
\subsection{Equal multitudes}\label{determining}

Only in some cases is it  possible to determine that two sets have the same number of elements.
\begin{quote} An equality of these multiplicities can only be concluded if some other reason is added, such as that both multitudes have exactly the same \emph{determining ground, [Bestimmungsr\H unde]}, 
e.g. they have exactly the same \emph{way of being formed [Entstehungsweise]}. (PU \S 21). \end{quote} 

Bolzano does not explain the exact meaning of the term \emph{determining ground} here. Its general mathematical definition is not precisely established (Bellomo \& Massas 2021, p. 12).  From other Bolzano's texts it follows that to \emph{determine an object} means to describe all representations that the object falls under. The determination is complete if the representation of an object is unique (\v{S}ebest\'{i}k 1992, p. 460). 

Bolzano provides two examples of equal infinite multitudes that have the same determining ground. The first example is the multitude of quantities which lie between two given quantities having the same distance, see (PU \S 29).\footnote{The word quantities here means number quantities, i.e. rational or measurable numbers that are actually very similar to real numbers. (Russ \& Trlifajov\' a 2016), (Fuentes Guill\' en 2021).} 
\begin{quote} It is no less clear that it will be found that the whole \emph{multitude} of quantities which lie between two given quantities, e.g. $7$ and $8$, \dots depends solely on the magnitude of the distance of those two boundary quantities from one another, i.e. on the quantity $8-7$, and therefore must be equal whenever this distance is equal. 
\end{quote}
 The second example is geometrical. Bolzano claims that every straight line, which is \enquote{not only similar to another but also \emph{geometrically} equal (i.e. coincides with it in all \emph{characteristics which are conceptually representable} through comparison with a given distance) must also have an equal multitude of points}. (PU \S 49.1). He even goes so far as to express an exact multitude of points relative to the unit line. He does not forget the endpoints; where they overlap, he subtracts them. 
\begin{quote} If we designate the multitude of points that lie between $a$ and $b$, including $a$ and $b$, by $E$, and take the straight line $ab$ as the unit of all \emph{lengths} then the multitude of points in the straight line $ac$, which has the length $n$ (a whole number) if its endpoints $a$ and $c$ is to be counted in, is $= nE - (n-1)$. (PU \S 49.3). \end{quote} 


According to Bolzano, two intervals of the same length have the same multitude of quantities. Likewise, two straight lines of the same length have the same multitude of points. 
This is perfectly consistent with PW. 

\subsection{Infinite Series}

Bolzano had already investigated  \emph{infinite quantity expressions} in the 7\textsuperscript{th} section of his \emph{Pure Theory of Numbers}, which was written in the early 1830s 
where he introduced \emph{measurable numbers} on their basis (Russ \& Trlifajov\'{a} 2016, p. 46). While here he mainly worked with infinite expressions which can be interpreted as Cauchy sequences; later on, in \emph{Paradoxes of the Infinite} he dealt particularly with infinite series which can be interpreted as divergent sequences. The following examples are from \S 29 and \S 33 of PU.\footnote{For the sake of clarity, we use the designation $N$ instead of Bolzano's $\overset{0}{N}$, $M$ instead of $\overset{m}{N}$, $P$ instead of $\overset{1}{S}$, and $S$ instead of $\overset{2}{S}$.} 
\begin{itemize}
\item $N\ = 1 + 1 + 1 + 1 + \dots  \text{in inf.} $  
\item $M\ = \underbrace{ \dots  }_m  1 + 1 + 1 + \dots  \text{in inf.} $, where the first $m$ terms are omitted. 
\item $P \ = 1 + 2 + 3 + \ 4 + \dots  \text{in inf.} $
\item $S \ = 1 + 4 + 9 + 16 + \dots  \text{in inf.} $

\end{itemize}

While in the usual mathematical understanding all these series have the same infinite sum, for Bolzano they are examples of various infinite quantities. Several important principles for adding, subtracting, and comparing them follow from the way Bolzano treats and comments on them.

\begin{enumerate}[(B1)]

\item All infinite series have one and the same multitude of terms.\footnote{\enquote{The multitude of terms in both series [$B$ and $C$] is certainly  \emph{the same}. By the fact that we raise every single term of the series $B$ to the square into the series $C$, we alter merely the nature (the magnitude) of this terms and not their plurality.}
(PU, \S 33).} 
\item Unless explicitly stated otherwise, as is the case of $M$. Then the difference between $M$ and $N$ is just $m$,   
$N - M = m.\footnote{\enquote{If we remove $n$ terms from the series, then its sum is smaller by exactly the sum of the removed terms. Then we obtain by subtraction the certain and quite unobjectionable equation $\underbrace{1 + \dots +1}_m = m = N - M$.}
(PU, \S 29).}$
\item $S > P$ since every corresponding term of $S$ is greater than in $P$.\footnote{\enquote{But if the multitude of terms in $P$ and $S$ is the same, then it is clear that $S$ must be greater than $P$, since, with the exception of the \emph{first} term, each of the remaining terms in $S$ is definitely greater than the corresponding one in $P$.} (PU, \S 33).}  

\item  $S$ is infinitely greater, i.e. greater than every finite multiple of $P$, $S >> P$.\footnote{In \S 33, Bolzano says that if we subtract $P$ from $S$ we obtain \enquote{an infinite series with an equal numbers of terms as $P$, namely $$0, 2, 6, 12, 20, 30, 42, 56 \dots n(n-1) \dots  \text{in inf.}$$  in which, with the exception of the first term, all succeeding terms are greater than corresponding terms in $P$.  If we subtract from this remainder $P$ for the second times, then we obtain the second remainder \dots Now, since these arguments can be continued without the end it is clear that $S$  infinitely greater than $P$, while in general we have
$$S - P = (1-m)+(2^2-2m) + (3^3-3m) + \dots(m^2-m^2) + \dots +m(n-m) + \dots \text{in inf.}$$
In this series only a finite multitude of terms is negative, namely the first $m-1$ are negative and the $m$-th is $0$, but all succeeding ones are positive and increase indefinitely.}}  
\item If we change the order of finitely many terms of the series the quantity doesn't change.\footnote{Unlike other principles, this is not literally exact. In \S 32, Bolzano says: \enquote{In particular, a \emph{series}, if we want to consider it only as a quantity, namely as the \emph{sum} of its terms must \dots have such a nature that it undergoes no change in value when we make a change in order of its terms, it must be that $$(A+B)+C = A+(B+C)=(A+C)+B.\text{''}$$}} 

\end{enumerate}
According to Bolzano, the series $N$ represents  the multitude of all natural numbers and $M$ the multitude of natural numbers greater than $m$ (PU, \S 29).

\subsection{Interpretation}
In two recent works (Trlifajov\' a, 2018), (Bellomo \& Massas 2021), Bolzano's infinite series are interpreted as sequences of partial sums to demonstrate that they form a consistent system.\footnote{Divergent series are treated similarly in (Bartlett, Logan \& Nemati 2020).} 
By converting the series to sequences, the requirement (B1) of \enquote{one and the same multitude} of terms is satisfied. Each term of a sequence has its specific uniquely determined place. Two sequences are compared, summarized and subtracted \emph{componentwise}, i.e. according their corresponding terms, which is in agreement with the justification of (B3) and (B4). The omitted terms of series, as in the case of $M$, are simply replaced by zero, so specific differences, as between $N$ and $M$, are naturally preserved, see (B2). 
\begin{definition}
Let $a_i \in \mathbb N_0$ for $i \in \mathbb N$.\footnote{As usually we denote natural numbers $\mathbb N = \{1, 2, 3, \dots \}$ while $\mathbb N_0 = \{0,1,2,3, \dots \}$ includes zero.} We interpret Bolzano's series $$a_1 + a_2 + a_3 + \dots \text{ in inf.}$$ 
as the sequence $(s_n)_{n \in \mathbb N}$ where $s_n = a_1 + \dots + a_n$ for any $n \in \mathbb N$
$$a_1 + a_2 + a_3 + \dots \text{ in inf.} \sim (s_1, s_2, s_3, \dots ) = (s_n)_{n \in \mathbb N}.$$
\end{definition}
Let $\mathbb N_0^\mathbb N$ be the set of sequences of natural numbers. If there is no danger of an error we will designate a sequence $(a_n)_{n \in \mathbb N} \in \mathbb N_0^\mathbb N$ just as $ (a_n)_n$ or only $(a_n)$.  

In this interpretation
\begin{itemize}
\item $N = 1 + 1 + 1 + 1 + \dots \text{ in inf.} \sim (1, 2, 3, \dots) = (n)_n$
\item $M = \underbrace{0 +\dots + 0}_m + 1 + 1 + \dots \text{ in inf.} \sim ( \underbrace{0, \dots, 0}_m, 1, 2, 3, \dots ) = (\max \{0, n-m\})_n$
\item $P = 1 + 2 + 3 + 4 + \dots \text{ in inf.} \sim (1, 3, 6, 10, \dots ) = (\frac {n \cdot (n+1)}{2})_n$
\item $S = 1 + 4 + 9 + 16 + \dots \text{ in inf.} \sim (1, 5, 16, 32, \dots ) = (\frac{n(n+1)(2n+1)}{6})_n$

\item $m = \underbrace{1 + \dots +1}_m + 0 + 0 +\dots \text{ in inf.} \sim (1, 2, 3, \dots, m, m, m, \dots)$ for $m \in \mathbb N$.  
\end{itemize}
Any Bolzano's series of natural numbers is now uniquely interpreted as a non-decreasing  sequence of natural numbers and vice versa, any non-decreasing sequence can be written as Bolzano's series. The sum and the product of two sequences is defined componentwise in accordance with adding and subtracting series, see (B3) and the note on (B4).\footnote{There is no direct example of multiplication of infinite series in PU. We assume that product as well as sum and difference should be the extension of the product of a finite number of terms.
$$(a_1 + a_2 + \dots \text{ in inf.}) + (b_1 + b_2 + \dots \text{ in inf.}) \sim (\sum_{i=1}^{n}a_i)_n + (\sum_{i=1}^{n}b_i)_n = (\sum_{i=1}^{n}a_i + \sum_{i=1}^{n}b_i)_n =  (\sum_{i=1}^{n}(a_i + b_i))_n$$
$$(a_1 + a_2 + \dots \text{ in inf.}) \cdot (b_1 + b_2 + \dots \text{ in inf.}) \sim (\sum_{i=1}^{n}a_i)_n \cdot (\sum_{i=1}^{n}b_i)_n = (\sum_{i=1}^{n}a_i \cdot \sum_{j=1}^{n}b_i)_n =  (\sum_{i,j =1}^{n}(a_i \cdot b_j))_n$$}
\begin{definition}\label{Bdef}
 Let $(a_n), (b_n)$ be two sequences of natural numbers. 
$$(a_n) + (b_n) = (a_n + b_n).$$
$$(a_n) \cdot (b_n) = (a_n \cdot b_n).$$
\end{definition}

According to (B5), the value of a series undergoes no change when we change the order of finitely many terms, which means that two series are equal if terms of sequences which are their interpretations are equal sarting from a sufficiently great index. In this sense, we define an \emph{equality} $=_\mathcal F$ and an \emph{ordering} $ <_\mathcal F$ of sequences.\footnote{The justification of the designation $=_\mathcal F$ with the index $\mathcal F$ will be in Theorem \ref{just}.} We also define the relation to be \emph{infinitely greater} $ >>_\mathcal F$ as being greater than every finite multiple. And finally, two sequences are \emph{equal in order}, $\approx_\mathcal F$, if for each of them its finite multiple is greater than the other sequence.

\begin{definition}\label{Bdef} Let $(a_n), (b_n)$ be two sequences of natural numbers, $k,m,n \in \mathbb N$. 
$$(a_n) =_\mathcal F (b_n) \text{ if and only if } (\exists m)(\forall n)(n > m \Rightarrow a_n = b_n).$$
$$(a_n) <_\mathcal F (b_n) \text{ if and only if }(\exists m)(\forall n)(n > m \Rightarrow a_n < b_n).$$
$$(a_n) <<_\mathcal F (b_n) \text{ if and only if }(\forall k)(k \cdot (a_n) <_\mathcal F (b_n)).$$
$$(a_n) \approx_\mathcal F (b_n) \text{ if and only if } (\exists k)(k \cdot (a_n) \geq_ \mathcal F (b_n)) \wedge (\exists k)(k \cdot (b_n) \geq_\mathcal F (a_n)).\footnote{In Bachmann-Landau notation, $(a_n) <<_\mathcal F (b_n)$ has the same meaning as $(a_n) = o(b_n)$ and $(a_n) \approx_\mathcal F (b_n)$ as $(a_n) = \Theta(b_n)$. This  relational notation is old. Paul du Bois-Reymond already had used it in \emph{Sur la grandeur relative des infinis des fonctions} in 1871 (Knuth 1976, p. 21).}$$
Equality $=$ of two sequences without the subscript means that all terms of the sequences are equal. It is the same for ordering $<$.
\end{definition}

\begin{theorem}\label{just} The result of the equality $ =_\mathcal F$ and order $ <_\mathcal F$ is the same as if we define equality and order by the Fr\'{e}chet filter on natural numbers. Both  $=_\mathcal F$ and $\approx_\mathcal F$ are equivalence relations.  The relations $<_\mathcal F$ and $<<_\mathcal F$ are strict partial orders. 
\begin{proof} The Fr\'{e}chet filter is the set of all complements of finite subsets of natural numbers. $$\mathcal F = \{A \subseteq \mathbb N; \mathbb N \setminus A \textrm{ is finite}\}$$ Let $(a_n), (b_n)$ be two sequences of natural numbers. Their equality and order by the Fr\'{e}chet filter is defined $$(a_n) =_\mathcal F (b_n) \text{ if and only if } \{n; a_n = b_n\} \in \mathcal F.$$
 $$(a_n)  <_\mathcal F (b_n) \text{ if and only if } \{n; a_n < b_n\} \in \mathcal F.$$
 If two sequences are equal from a sufficiently large index, then they differ only in finitely many terms. Thus, they are equal by the  Fr\'{e}chet filter. Vice versa, if two sequences are equal by the Fr\'{e}chet filter then they must be equal from a sufficiently large index. The same applies for the order. The symbols $=_\mathcal F$ and $<_\mathcal F$ are hereby justified. 

Symmetry and reflexivity of $=_\mathcal F$ and $\approx_\mathcal F$ are evident from Definition \ref{Bdef}, transitivity is easily proven.  
Irreflexivity and transitivity of  $<_\mathcal F$ and $<<_\mathcal F$ are also obvious from Definition \ref{Bdef}. 
\end{proof}
\end{theorem}

\begin{theorem}\label{ring}
Let $S = \{(a_n), a_n \in \mathbb N \wedge (\forall n)(a_n \leq a_{n+1})\}$ be the set of non-decreasing sequences of natural numbers. Then the structure $(S, + , \cdot, =_\mathcal F, <_\mathcal F)$ where the equality and the ordering is defined by Definition \ref{Bdef} is a partial ordered non-Archimedean commutative semiring.\footnote{Moreover, we can define  $$(a_n) \leq_\mathcal F (b_n) \text{ if and only if }(\exists m)(\forall n)(n > m \Rightarrow a_n \leq b_n).$$ Then the structure $(S, \leq_\mathcal F)$ is a distributive lattice. The greatest lower bound and lowest upper bound of $(a_n)$ and $(b_n)$ in $S$ are found by taking component-wise greatest lower bound and lowest upper bound in the underlying lattice $(\mathbb N, \leq)$. Thank you Ansten M\o{}rch Klev for this proof. }

\begin{proof}
The sum and the product of non-decreasing sequences of natural numbers are defined componentwise. Hence the properties of associativity, commutativity, distributivity and the existence of the neutral element $(0)_n$ and the unit $(1)_n$ are valid both for addition and multiplication. Partial ordering follows from Theorem \ref{ring}. The sequence  $(n^2)_n$ 
is infinitely greater than the sequence $(n)_n$ 
that is infinitely greater than any constant sequence.  
\end{proof}
\end{theorem}
Under this interpretation, Bolzano's infinite quantities form a consistent and meaningful algebraic structure and Bolzano's assertions concerning infinite series  are valid. (Trlifajov\' a 2018, p. 698).\footnote{Compared to this interpretation, Bellomo and Massas argue that the product of two sequences should not be interpreted componentwise. (Bellomo \& Massas 2021, p. 36). The main reason is based on the passage in wich Bolzano presents an example of an infinite quantity of a \emph{higher order} 
 $$N + N  + N + \dots \textrm{ in inf.} = N^2$$
The interpretation of $N$ is $(1,2,3, \dots)$ and that of $N^2$ would be $(1,4,9, ... ) = (n^2)$ in componentwise interpretation. But the same sequence is the interpretation of the series of odd numbers $1 + 3 + 5 + \dots \textrm{ in inf.}$ According to authors, this is impossible. $N$ is surely greater than any odd numbers, so by (B3) $N^2$ has to be greater than the sum of odd numbers, and not equal. However, 
the series $N + N  + N + \dots \textrm{ in inf.}$ is just symbolical, Bolzano does not admit series of infinite quantities, he says that \enquote{it is absurd to speak of the last term of the series that has the value $N$}. Our explanation why $N^2$ is equal to $1 + 3 + 5 + \dots \textrm{ in inf.}$ will be given in Definition \ref{def} and  Example \ref{ex1}, 5. 


However, though the authors present a slightly different interpretation based on an iterative ultraproduct they show that Bolzano's theory of infinite quantities can be consistently interpreted as a non-commutative ordered ring. (Bellomo \& Massas 2021, p. 42).
}

Bolzano's aim was primarily to develop an arithmetic of infinite series of integers (Bellomo \& Massas 2021, p. 49).  I do not claim that Bolzano intended to use them to determine sizes of infinite countable sets, however, since his theory preserves PW we will use them and Bolzano's principles for this purpose.

\section{Sizes of Canonically Countable Sets}\label{size}

\subsection{Concept of Size}

If we wish to count a great number of things, or even an infinite number, it may be convenient to arrange them in smaller groups according to some rule, count the number of their elements, and then add these numbers. 
A nice popular example of Numerosity Theory (NT) is the world population. If one wants to count them, it is convenient to divide them by nation, count the population of each nation, and finally add up the numbers of nations, alphabetically for example. 
We must have a criterium for how to make this arrangement.

\begin{quote} If we want to apply this idea to infinite sets, it is necessary to have a criterium to arrange elements into smaller groups. (Benci \& Di Nasso 2019, p. 277) \end{quote} 

We use the same concept for countable sets. 
We arrange a set $A$ into a sequence of mutually disjoint finite subsets of $A$ 
such that their union is the whole set. 
\begin{definition} A sequence $(A_n)_{n \in \mathbb N}$ of finite disjoint subsets of $A$, $A_n \subseteq A$, is an \emph{arrangement of a set $A$} if
$$A = \bigcup\{A_n; n \in \mathbb N\}$$ 
The sets $A_n$ are called \emph{components} of the arrangement. 
\end{definition} 
Of course, there are infinitely many, even uncountably many, possible arrangements of a countable set into a sequence of finite components. But we shall determine a \emph{canonical} arrangement in the following Section \ref{ca}. This will depend on the \emph{determining ground} of a set. A size of a canonically arranged set $A$ will be defined as the sum of finite cardinalities of its components, Bolzano's infinite series. 
$$|A_1| + |A_2| + |A_3| + \dots \text{ in inf.}$$
This is interpreted as a sequence of partial sums, a \emph{size sequence} of $A$.
$$\sigma(A) = (|A_1|, |A_1| + |A_2|, |A_1| + |A_2| + |A_3|, \dots). $$ 
\begin{remark}\label{ord}
\begin{enumerate}
\item Canonically arranged sets are \emph{countable} by Cantor's definition. But it does not have to work the other way around. In particular, this applies to ordinal numbers. These notions are based on Hume's Principle, which is in contradiction with PW. We will not deal here with their sizes. The same ordinal number $\omega$ has the set of natural numbers and all its infinite subsets. Ordinal numbers are rather types of well-orderings; it is meaningless to determine their sizes. 
\item Our method is close to that of Numerosity Theory (NT). 
We use slightly different terminology on purpose to point out some different assumptions of our theory. \emph{Canonically countable} sets are essentially identical to their \emph{labelled} sets.\footnote{Label of an element of a set determines the index of a component} (Benci \& Di Nasso 2019, p. 278). More in Section \ref{num}. 

\end{enumerate}
\end{remark}

\subsection{Canonical Arrangement}\label{ca}

Our aim is to describe a way how to unambiguously arrange a countable set into a sequence of finite components according to its  \emph{determining ground}. This method should lead to the correct result even when applied to a finite set. 

The natural ordering of the natural numbers $\mathbb N$ is almost literally based on their \emph{determining ground}. Its canonical arrangement is defined so that the $n$-th component contains just the number $n$, it is a one-element set $\{n\}$. This arrangement forms a kind of skeleton which will be gradually extended according to certain rules.

First, subsets and supersets of canonically arranged sets should preserve the same arrangement. If $A, B$ are canonically arranged sets, $(A_n)_{n}$, $(B_n)_{n}$ are their canonical arrangements and $B \subseteq A$ then it should also be that $B_n \subseteq A_n$ for all $n$. 
Consequently, the canonical arrangement of a subset of a canonically arranged set is given uniquely. Moreover, a canonical arrangement of supersets can only be defined with regard to this rule.

The second rule concerns a canonical arrangement of the Cartesian product $A \times B$ of two canonically countable sets $A$ and $B$. 
How to arrange it into finite disjoint subsets so that the same method also leads to the correct result in the case of finite sets? The solution is a good old way but with a new meaning. The $n$-th component contains a union of $A_i \times B_j$ where $n = \max\{i,j\}$, the \enquote{n-th frame} in this picture. 
\begin{center}
\catcode`\-=12
\begin{tabular}{c||c c c c c c c c c c}
& $A_1$ & $A_2$ & $A_3$ & $A_4$ & $\dots$\\
\hline
\hline
$B_1$ & \ $A_1\times B_1$ \.& \ $A_2\times B_1$ \. & \ $A_3\times B_1$ \. & \ $A_4\times B_1$ \. & \dots \\
\cline{2-2}
$B_2$ & $A_1\times B_2$ & \ $A_2\times B_2$ \.& \ $A_3\times B_2$ \. & \ $A_4\times B_2$ \. & \dots \\
\cline{2-3}
$B_3$ & $A_1\times B_3$ & $A_2\times B_3$ & \ $A_3\times B_3$ \. & \ $A_4\times B_3$ \. & \dots \\
\cline{2-4}
$B_4$ & $A_1\times B_4$ & $A_2\times B_4$ & $A_3\times B_4$ & $A_4\times B_4$ \ \. & \dots\\
\cline{2-5}
$\dots$ &\dots &\dots & \dots &\dots & \dots \\
\end{tabular}
\end{center}
\begin{itemize}
\item It is an arrangement of $A \times B$ into finite components that covers it all.\footnote{The same way of arrangement of the product of divergent series is proposed and  justified in the paper \emph{The Products of Hyperreal Series and the Limitations of Cauchy Products} (Bartlett 2022).}  
$$A \times B = \bigcup\{(A \times B)_n; n \in \mathbb N\}.$$ 
\item It can be considered as an extension of the arrangements of both $A$ and $B$, so it satisfies the previous rule. In fact, if we assume that $A$ has the same arrangement as $A \times \{1\}$ then
$$A_n \times \{1\} = (A \times \{1\})_n \subseteq (A \times B)_n.$$ 
\item The same method applied to finite sets gives the correct result no matter what their arrangement is.\footnote{A different arrangement is usually proposed, such that $$(A \times B)_n = \bigcup\{A_i \times B_j, n = i+j\}.$$
But in this case, it would only cover a half of the Cartesian product and not the whole.}
\end{itemize}

\begin{definition}\label{def} A canonical arrangement is defined successively. Let $n \in \mathbb N$. 
\begin{enumerate} 
\item The $n$-th component of canonical arrangement of natural numbers $\mathbb N$ is 
$$\{n\}.$$
\item Let $(A_n)_n$ be the canonical arrangements of $A$ and $B \subseteq A$. Then $B$ is also canonically arranged and its $n$-th component 
$$B_n = A_n \cap B.$$ 
\item 
 Let $(A_n)_n$, $(B_n)_n$ be the canonical arrangements of $A$ and $B$. The $n$-th component of the canonical arrangement of the Cartesian product $A \times B$ is 
$$(A \times B)_n = \bigcup\{A_i \times B_j, n = \max\{i,j\}\}.$$

\end{enumerate}
Canonically arranged sets will be also called \emph{canonically countable}.
\end{definition}

\subsection{Characteristic and Size Sequences}
A canonical arrangement of a set is determined by a sequence of its subsets. We introduce two notions, characteristic and size sequence, that describe this arrangement. The  $n$-th term of the \emph{characteristic sequence} is the cardinality of its $n$-th component, it can be zero. The \emph{size sequence} is a sequence of partial sums of the characteristic sequence; this expresses the \emph{size} of a set, its set-size.

\begin{definition} Let $A$ be a canonically countable set and $(A_n)_n$ its arrangement.  
\begin{itemize} 
\item A \emph{characteristic sequence} of $A$ is the sequence $\chi(A) = (\chi_n(A))_n \in \mathbb N_0^\mathbb N$ where $$\chi_n(A) = |A_n|.$$
\item A \emph{size sequence} of $A$
is the sequence $\sigma(A) = (\sigma_n(A)) \in \mathbb N_0^\mathbb N$ where $$\sigma_n(A) =  \chi_1(A) + \dots + \chi_n(A).$$
\end{itemize}
\end{definition}

\begin{remark}
A simple relationship follows directly from the definition. For a set $A$ 
$$\sigma_n(A) = \sigma_{n-1}(A) + \chi_n(A).$$ 
\end{remark}

\begin{example}\label{ex1}

Characteristic and size sequences.
\begin{enumerate}
\item Natural numbers $\mathbb N = \{1, 2, 3, \dots\}$
$$\chi(\mathbb N) = (1,1,1,1, \dots), \quad \sigma(\mathbb N) =  (1,2,3,4, \dots) .$$
\item If $A \subseteq \mathbb N$ then $\chi_n(A)$ is equal either $1$ or $0$. 
$$\chi_n(A) = 1 \Leftrightarrow n \in A.$$
\item The two-elements subsets $\{1, 2\}$ and  $\{3, 4\}$ of $\mathbb N$:
$$\chi( \{1, 2\}) = (1,1,0,0,0,0, \dots), \quad\sigma( \{1, 2\}) = (1,2,2,2,2,2, \dots).$$ 
$$\chi( \{3, 4\}) = (0,0,1,1,0,0 \dots), \quad \sigma( \{3, 4\}) = (0,0,1,2,2,2, \dots).$$
\item Their complements: 
$$\chi(\mathbb N \setminus \{1, 2\}) = (0, 0, 1, 1, 1, \dots), \quad \sigma(\mathbb N \setminus \{1, 2\}) = (0, 0, 1, 2, 3, 4, 5\dots).$$
$$\chi(\mathbb N \setminus \{3, 4\}) = (1, 1, 0, 0, 1, \dots), \quad \sigma(\mathbb N \setminus \{3, 4\}) = (1, 2, 2, 2, 3, 4, 5\dots).$$

\item Even numbers $\mathbb E = \{2,4,6,8, \dots\}$ and odd numbers: $\mathbb O = \{1,3,5,7, \dots \}$. 
$$\chi(\mathbb E) = (0,1,0,1,0,1, \dots), \quad \sigma(\mathbb E) = (0,1,1,2,2,3,3 \dots).$$ 
$$\chi(\mathbb O) = (1,0,1,0,1,0 \dots), \quad \sigma (\mathbb O) = (1,1,2,2,3,3, \dots).$$ 

\item Cartesian product $\mathbb N \times \mathbb N$: $$(\mathbb N \times \mathbb N)_1 = \{[1,1]\}, (\mathbb N \times \mathbb N)_2 = \{[1,2], [2,2],[2,1]\}$$ 
$$(\mathbb N \times \mathbb N)_3 = \{[1,3], [2,3],[3,3],[3,2],[3,1]\}, \dots $$
$$\chi(\mathbb N \times \mathbb N) = (1, 3, 5, 7,  \dots) = (2n-1)_n, \   
\sigma(\mathbb N \times \mathbb N) = (1,4, 9, 16, \dots) = (n^2)_n.$$
\end{enumerate}

\end{example}

 \begin{definition}
The sequence $\sigma (\mathbb N)$ represents the size of natural numbers. We denote it symbolically as
$$\alpha = \sigma (\mathbb N)= (n)_n = (1, 2, 3, \dots)\footnote{NT uses the same symbol; $\alpha$ is defined as a \enquote{new} number which can be considered as the sequence $(n)_n$ (Benci \& Di Nasso 2019, p. 15).}$$ 
A constant sequence $(k)_n$ is identified with the number $k \in \mathbb N$ as usual. 
$$(k)_n = k$$
\end{definition}

\begin{remark}
The sequence $\alpha$ is the interpretation of Bolzano's series 
$$N = 1 + 1 + 1 + \dots \text{ in inf.} \sim \alpha$$
that represents the \enquote{multitude} of all natural numbers. However, there is a significant difference between this $\alpha$ and Cantor's $\omega$ that is the smallest infinite number. The infinite quantity $\alpha$ is neither the limit of finite numbers nor the smallest infinite number. It can be added, subtracted, multiplied while preserving laws of arithmetic. 
\end{remark}
\subsection{Structure of Size Sequences}

Size sequences are interpretations of Bolzano's series. 
By Definition \ref{Bdef} the sum and the product of two sequences are defined componentwise. The equality $=_\mathcal F$ and order $ <_\mathcal F$ are defined modulo the Fr\' echet filter, that in practice means that these relations are valid from a sufficiently great index. 
Similarly, the relations of infinite order $<<_\mathcal F$ and of equality in order $\approx_\mathcal F$ are defined modulo the Fr\' echet filter. 
According to Theorem \ref{ring} this structure is a \emph{partial ordered non-Archimedean semiring.}  

\begin{example}  We use the same notation as in Example \ref{ex1}. 
\begin{itemize}

\item $\sigma(\{1, 2\}) =_\mathcal F  \sigma(\{3, 4\}) =_\mathcal F 2$; \quad $\sigma(\mathbb N \setminus \{1, 2\}) + \sigma(\{1, 2\})  =  \alpha $
\item $\alpha <_\mathcal F  \alpha^2$  and even $\alpha <<_\mathcal F  \alpha^2$ for $(\forall k)(k\cdot \alpha <_\mathcal F \alpha^2)$.
\item $\sigma(\{1, 2\}) <<_\mathcal F \alpha$  since  $(\forall k \in \mathbb N)( k\cdot \sigma(\{1, 2\}) =_\mathcal F k \cdot 2 <_\mathcal F \alpha)$   

\item $\sigma(\mathbb O) + \sigma(\mathbb E) = \alpha; \quad \sigma(\mathbb E) \leq \sigma(\mathbb O) \leq \sigma(\mathbb E) + 1$; \quad $\sigma(\mathbb E) \approx_\mathcal F \sigma(\mathbb O) \approx_\mathcal F \alpha$;

\end{itemize}
\end{example}

Every element of a canonically countable set has a unique place in a clearly defined component. If two such sets $A, B$ have a non-empty intersection then for all its elements $x \in A \cap B$ it holds that $x \in A_n \Leftrightarrow x \in B_n$. From this uniqueness follows:

\begin{theorem}\label{union}
Let $A = \bigcup\{A_n, n \in \mathbb N\}, B = \bigcup\{B_n, n \in \mathbb N\}$ be two canonically countable sets. Then the size of their union is
$$\sigma(A \cup B) = \sigma(A) + \sigma(B) - \sigma(A \cap B)$$
\end{theorem}
\begin{proof}
The $n$-th component of the arrangement is $(A \cup B)_n = A_n \cup B_n.$ Consequently 
$$\chi_n(A \cup B) = |A_n \cup B_n| = |A_n| + |B_n| - |A_n \cap B_n|.$$ 
$$\sigma_n(A \cup B) = \Sigma_{i =1}^n \chi_i(A_i \cup B_i) = \Sigma_{i =1}^n(|A_i| + |B_i| - |A_i \cap B_i|) = \sigma_n(A) + \sigma_n(B) - \sigma_n(A \cap B).$$ 
\end{proof}
The following theorem shows that the size of the Cartesian product of two sets corresponds to the product of their size sequences. 
\begin{theorem}\label{union}
Let $A = \bigcup\{A_n, n \in \mathbb N\}, B = \bigcup\{B_n, n \in \mathbb N\}$ be two canonically countable sets. Then the size of their Cartesian product is
$$\sigma(A \times B) = \sigma(A) \cdot \sigma(B).$$\end{theorem}
\begin{proof}
The $n$-th component of their product $$(A \times B)_n = \bigcup\{A_i \times B_j, n = \max\{i,j\}\}.$$
The $n$-th term of the characteristic sequence $$\chi_n(A \times B) = \sum\{|A_i| \cdot |B_j|, \max\{i,j\} = n\}.$$ 
The $n$-th term of the size sequence $$\sigma_n(A \times B) = \chi_1(A \times B) + \dots + \chi_n(A \times B) = \sum\{|A_i| \cdot |B_j|, \max\{i,j\} \leq n\}= \sigma_n(A) \cdot \sigma_n(B).$$

\end{proof}

\begin{example}  By Example \ref{ex1}, the sizes of even and odd numbers are
$$\sigma (\mathbb E) = (0,1,1,2,2,3 \dots), \sigma (\mathbb O) = (1,1,2,2,3,3 \dots).$$ 
\begin{itemize}

\item $\sigma(\mathbb E \times \mathbb O) =\sigma (\mathbb E) \cdot \sigma (\mathbb O) = \sigma(\mathbb O \times \mathbb E) = (0, 1, 2, 4, 6, 9 \dots)$
\item $\sigma(\mathbb E \times \mathbb E) = \sigma (\mathbb E)^2 = (0, 1, 1, 4, 4, 9 \dots)$
\item $\sigma(\mathbb O \times \mathbb O) = \sigma (\mathbb O)^2 = (1, 1, 4, 4, 9, 9,  \dots)$ 
\item $\sigma(\mathbb N \times \mathbb N) = 2 \cdot \sigma(\mathbb E \times \mathbb O) + \sigma(\mathbb E \times \mathbb E) + \sigma(\mathbb O \times \mathbb O) = (1, 4, 9, 16, 25 \dots) = \alpha^2.$

\end{itemize}
\end{example}
The arrangement of a \emph{finite set} is not important. If it is a subset of a canonically countable set then its size sequence is equal to its cardinality from a sufficiently great index. 
If it is not, we will assume without loss of generality, unless explicitly stated otherwise, that all its elements lie in the first component. 

\begin{definition}\label{remark}
If $A$ is a finite set that is not canonically arranged then we define it $$A_1 = A.$$ 
\end{definition}
In any case, the size of a finite set is equal to its cardinality $\sigma(A)=_\mathcal F |A|.$

\begin{theorem} Let $A, B$ be two canonically countable sets. 
\begin{enumerate}
\item 
Size sequences of finite sets are equal to their cardinalities. If $A$ is finite then $$\sigma(A) =_\mathcal F |A|.$$
\item 
The Part-Whole Principle is valid. If $A$ is a proper subset of $B$, $A \subset B$, then $$\sigma(A) < \sigma(B).$$
\item The ordering is \emph{discrete} that means
$$\sigma(A) < \sigma(B) \Rightarrow \sigma(A \cup \{x\}) \leq \sigma(B).\footnote{This says that two sets cannot differ in size by less than a whole element (Parker 2013, p. 594).}$$
\end{enumerate}
\end{theorem}
\begin{proof}
\begin{enumerate}
\item See previous Definition \ref{remark}. 

\item If $A \subset B$ then $B \setminus A \neq \emptyset$ and $\sigma(B) = \sigma(A) + \sigma(B \setminus A)$ according to Theorem \ref{union}. 
Consequently $\sigma(A) < \sigma(B)$. 

\item Discreteness follows from the fact that the terms of size sequences are natural numbers. If $\sigma_n(A) < \sigma_n(B)$ from a sufficienly great index then from the same index $\sigma_n(A\cup \{x\}) \leq \sigma_n(A) + 1 \leq \sigma_n(B)$.

\end{enumerate}
\end{proof}

\subsection{Partial or Linear Order?}\label{partial order}

The ordering of size sequences is only partial and not linear.\footnote{The ordering $<$ defined on a set $S$ is \emph{linear (total)} if for every $x,y \in S$ one of the three options holds true $x < y \vee x = y \vee y < x.$} 

\begin{example} For the sizes of even numbers $\sigma(\mathbb E)$ 
and odd numbers $ \sigma(\mathbb O)$, it holds 
$$\sigma(\mathbb E) \leq \sigma (\mathbb O) \leq \sigma(\mathbb E) + 1\footnote{The reason is that the first natural numbers is $1$, if it was $0$ it would be the other way around.} \quad \textrm{and} \quad \sigma(\mathbb E) + \sigma(\mathbb O) = \alpha.$$ 
It is not generally determined whether $\sigma(\mathbb O) = \sigma (\mathbb E)$ or $ \sigma(\mathbb O) = \sigma (\mathbb E) +1$.  If $\sigma(\mathbb O) = \sigma (\mathbb E)$ then $\alpha$ would be even, otherwise $\alpha$ would be odd. 
Thus, it is also not determined whether $\alpha$ is odd or even. And, of course, many such  things. 

 \end{example}

The reason is the use of the Fr\'{e}chet filter in Definition \ref{Bdef} of equality and order of sequences. If we used a non-principal ultrafilter instead as Numerosity Theory does, we would obtain a linear ordering. But the assignment of sizes would depend on the choice of an ultrafilter $\mathcal U$ that is always to some extent arbitrary (Benci \& Di Nasso 2003, p. 12), (B\l{}asczyk 2021, p. 88).

\begin{itemize}
\item If $\mathbb E \in \mathcal U$ then $\sigma(\mathbb O) =_\mathcal U \sigma(\mathbb E)$, thus $\alpha = 2 \cdot \sigma(\mathbb E)$ and $\alpha$ is even.
\item If $\mathbb E \notin \mathcal U$ then $\sigma(\mathbb O) =_\mathcal U \sigma(\mathbb E) + 1$, thus $\alpha = 2 \cdot \sigma(\mathbb E) + 1$ and $\alpha$ is odd. 
\end{itemize}

This problem can be overcome to some extent by using a selective ultrafilter containing predetermined sets. However, yet many of its properties remain arbitrary. We shall return to this question and discuss it in detail in Section \ref{num}. 

Our choice to use the Fr\'{e}chet filter leads to only a partial ordering. 
We admit that we cannot determine all relations and properties of size sequences. On the other side, we do not need ultrafilters or any other non-constructive mathematical objects including the Axiom of Choice. Our results are unambiguous and independent. Moreover, since the Fr\'{e}chet filter is contained in any non-principal ultrafilter everything we prove would be valid if we eventually modified Definition \ref{Bdef} and used the ultrafilter. 

\section{Number Structures}\label{numbers} 

\subsection{Natural Numbers and their Subsets}

Characteristic sequences of subsets of natural numbers contain only $0$ and $1$. Size sequences are non-decreasing sequences such that the following term is always either equal to the preceding one or to the preceding term plus 1. 
And vice versa, any such sequence uniquely determines a set of natural numbers having this size sequence.  

In ZF, subsets of natural numbers are often constructed by \emph{Separation Schema} using a function strictly increasing on the natural numbers. For example, square numbers are defined by the function $f(n) = n^2$.
$$\mathbb S = \{n \in \mathbb N; (\exists m)(m \in \mathbb N \wedge n = m^2)\}. $$
In that case the size sequence is a sequence of integer parts of the inverse function defined on real numbers.\footnote{The similar theorem about numerosities is Proposition 16.24  (Benci \& Di Nasso 2019, p. 284). } 

\begin{theorem}\label{previous} Let $f$ be a strictly increasing function defined on positive real numbers $\mathbb R^+ \longrightarrow \mathbb R^+$. Let $A \subseteq \mathbb N$ be such that 
$$A = \{n \in \mathbb N; (\exists m)(m \in \mathbb N \wedge n = f(m)\}.$$
Then for the $n$-th term of the size sequence of $A$ it is valid  
$$\sigma_n(A) = \floor{f^{-1}(n)}.\footnote{The integer part $\floor{x}$ of $x \in \mathbb R$ is defined as the greatest integer less than or equal to $x$  $$\floor{x} = \max\{n \in \mathbb Z; n \leq x\}.$$ 
}$$
\begin{proof}  If $f$ is a strictly increasing function defined on $\mathbb R^+$ then its inverse $f^{-1}$ is also a strictly increasing function defined on $\mathbb R^+$.  Its integer part $\floor{f^{-1}}$ is determined for all $n \in \mathbb N$ by the prescription 
$$\floor{f^{-1}(n)} = \max\{m \in \mathbb N; m \leq f^{-1}(n)\}.$$
 Then for any $n \in \mathbb N$
$$\sigma_n(A)  = |\{k \in \mathbb N; \ k \leq n \wedge (\exists m \in \mathbb N)(k = f(m)\}| = $$
$$ |\{m \in \mathbb N;\ m \leq f^{-1}(n)\}| = \max\{m \in \mathbb N; m \leq f^{-1}(n)\}= \floor{f^{-1}(n)}.$$

\end{proof}
\end{theorem}

\begin{theorem}
For any $p, k \in \mathbb N$ 
\begin{enumerate}
\item The size of all squares $\mathbb S$ is infinitely less than the size of all primes $\mathbb P$.\footnote{All the more, the size of the set of \emph{any} powers is infinitely smaller than the size of primes, since it is certainly smaller than the size of squares.  } 
\item The size of primes $\mathbb P$ is infinitely less than the size of all $k$-multiples $\mathbb M_k$ 

\item The size of all $k$-multiples $\mathbb M_k$ is equal in order to $\alpha$, $\sigma(M_k) \approx_\mathcal F  \alpha$.
\end{enumerate}
$$\sigma(\mathbb S_p) <<_\mathcal F  \sigma(\mathbb P) <<_\mathcal F \sigma(\mathbb M_k)\approx_\mathcal F \alpha$$

\begin{proof}

Let $\mathbb P = \{1, 2, 3, 5, 7,11, \dots \}$ be the set of primes. For any $n \in \mathbb N$, $\sigma_n(\mathbb P)$ is the number of primes less than or equal to $n$, it is usually denoted by $\pi(n)$ in the number theory.  
By the \emph{prime number theorem}, $\pi(n)$ is asymptotic to $\frac{n}{\log n}$ (Hardy \& Wright, p. 9). There are many estimates of $\pi(n)$ but this will do for now (Rosser \& Schoenfeld, 1962, p. 69). For $n\geq 17$\footnote{With the exception of $\pi(113)$.} 
$$\frac{n}{\log n} < \pi(n) < \frac{5 n}{4 \log n} =  1.25 \cdot\frac{ n}{\log n}$$
\begin{enumerate}
\item 

For each term of the size sequence of squares $\sigma(\mathbb S)$ it is valid
$$\sigma_n(\mathbb S) = \floor{\sqrt n}\leq \sqrt{n} < \frac{n}{\log n} < \pi(n)$$ 

Moreover, for any $m \in \mathbb N$ and from a sufficiently great $n$ 
$$m \cdot \sigma_n(\mathbb S) = m \cdot \floor{\sqrt n}\leq m \cdot \sqrt{n} < \frac{n}{\log n} < \pi(n)$$ 
Consequently, not only 
$\sigma(\mathbb S) <_\mathcal F \sigma(\mathbb P)$ \textrm{but also} $\sigma(\mathbb S) <<_\mathcal F \sigma(\mathbb P). $

\item Let $\mathbb M_k = \{k,2k,3k, \dots\}$ be the set of $k$-multiples. For its  $n$-th term $\sigma_n(\mathbb M_k)$ it is valid 
$$\sigma_n(\mathbb M_k) = \floor{\frac{n}{k}},  \quad \quad  \frac{n}{k}-1 \leq \sigma_n (\mathbb M_k) \leq \frac{n}{k}.$$
From a sufficiently great $n$ 
$$\pi(n) <  \frac{5n}{4\log n} < \frac{n}{k} -1< \sigma_n (\mathbb M_k)\quad \textrm{
for } \quad \frac{5}{4} \cdot k <  \log n - k \cdot \frac{ \log n}{n}. $$
Moreover, for any $m \in \mathbb N$ and from a sufficiently great $n$ 
$$m \cdot \pi(n) \leq m \cdot \frac{5 n}{4 \log n} < \frac{n}{k} -1< \sigma_n (\mathbb M_k)\quad \textrm{
for } \quad \frac{5}{4}\cdot mk < \log n - k \cdot \frac{ \log n}{n}. $$ 
Consequently, not only  
$\sigma(\mathbb P) <_\mathcal F \sigma(\mathbb M_k)$ but also $\sigma(\mathbb P) <<_\mathcal F \sigma(\mathbb M_k).$
\item From the inequality above follows that  
$n - \frac{1}{k} \leq k \cdot \sigma_n (\mathbb M_k) \leq n$, and so 
$$k \cdot \sigma_n (\mathbb M_k) \leq n \leq n - \frac{1}{k} + 1 \leq \sigma_n (\mathbb M_k) \cdot k + \sigma_n (\mathbb M_k) = (k+1) \cdot \sigma_n (\mathbb M_k).$$

\end{enumerate}
\end{proof}
\end{theorem}

\begin{theorem}\label{prop}
Let $\mathbb M_k$ be the sets of all $k$-multiples of natural numbers, $k \in \mathbb N$. Then 
$$ \alpha - k \leq k \cdot \sigma(M_k) \leq \alpha$$
\end{theorem}

\begin{proof} 
According to Theorem \ref{previous}, $\sigma_n(\mathbb M_k) = \floor{\frac{n}{k}}$ for any $n \in \mathbb N$. Thus $\frac{n}{k}-1 \leq \sigma_n (\mathbb M_k) \leq \frac{n}{k}.$ Consequently, $$n - k \leq k \cdot \sigma_n (\mathbb M_k) \leq n.$$ 
\end{proof}


\subsection{Integers}\label{integer}
\begin{theorem}
Both non-negative integers $\mathbb N_0$ and integers $\mathbb Z$ are canonically countable. Their sizes are
$$\sigma(\mathbb N_0) = \alpha + 1, \quad \sigma(\mathbb Z) = 2 \alpha + 1.$$ 
\begin{proof}
Since $\mathbb N_0 = \mathbb N \cup \{0\}$, so $\sigma(\mathbb N_0) = \sigma(\mathbb N) + \sigma(\{0\}) = \alpha + 1.$

Negative integers $\mathbb N^{-}$ can be represented as the Cartesian product $\mathbb N \times \{1\}$, $-n$ is represented as $(n,1)$. 
They have the same arrangement, hence the size sequence $$\sigma(\mathbb N^{-}) = \alpha.$$
Integers are the union of positive and negative integers, and zero $\mathbb Z = \mathbb N \cup \mathbb N^{-} \cup \{0\}$  
$$\sigma(\mathbb Z) = \sigma(\mathbb N) + \sigma(\mathbb N^{-}) + \sigma(\{0\}) = 2 \alpha + 1.$$
\end{proof}
\end{theorem}

\subsection{Rational Numbers}\label{rational}

We are looking for a representation of the rational numbers that respects their determining ground, extends the arrangement of the natural numbers, and in addition preserves the property required by Bolzano, see Section \ref{determining}, namely that two intervals of the same length have the same size. We call it \emph{homogeneity} of size. 

At the beginning of this Section, we introduce a canonical arrangement of rational numbers that satisfies the first two properties, and at the end, we prove its homogeneity.

\subsubsection{Unit Interval}
Let us start with a half-open unit interval of rational numbers  $\mathbb I = (0,1]_\mathbb Q$. Every $x \in \mathbb I$ can be uniquely expressed as a ratio  $x = \frac{k}{m}$ of two coprime natural numbers $k, m$ such that $k \leq m$.\footnote{We deal with the \emph{half-open} interval because we could simply put more of these side by side to form a larger interval.  And we deal with $(0,1]_\mathbb Q$ and not $[0,1)_\mathbb Q$, that would lead to the same result, just because the common use of the term \emph{coprime} numbers - 
i.e. their greatest common divisor is $1$.} 

When we divide an interval or a line segment we usually think of \enquote{halves} first, then \enquote{thirds}, \enquote{quarters}, etc. It is quite natural that the first component of a canonical arrangement of rational numbers from $\mathbb I$ contains $1$, because this is the only integer, the second component contains $\frac{1}{2}$, the third $\frac{1}{3}$, $\frac{2}{3}$, the fourth $\frac{1}{4}$, $\frac{3}{4}$, etc. The fourth component does not contain $\frac{2}{4}$ for it is already included in the second one. 

If we represent $\mathbb I$ as a subset of $\mathbb N \times \mathbb N$ of pairs of coprime natural numbers, such that the first term is less or equal than the second one we get exactly this arrangement. It is reasonable and corresponds to its determining ground.

\begin{definition} We introduce a predicate $\textrm{co}(k,m)$ saying that $k,m$ are coprime natural numbers.
$$\textrm{co}(k,m) \Leftrightarrow k,m \in \mathbb N \wedge k,m \textrm{ are coprime }.$$ The half-open unit interval of rational numbers  $\mathbb I = (0,1]_\mathbb Q$ is represented by the set
$$\mathbb I  = \{\frac{k}{m} \sim (k,m) \in \mathbb N \times \mathbb N; \textrm{co}(k,m) \wedge k \leq m\} \}.$$ 
\end{definition}

The Cartesian product $\mathbb N \times \mathbb N$ is canonically countable as well as its subset $\mathbb I \subseteq \mathbb N \times \mathbb N$. Its arrangement is determined 
by Definition \ref{def}.
$$\mathbb I_n = (\mathbb N \times \mathbb N)_n \cap \mathbb I = \{(k, m);n = \max\{k, m\} \wedge \textrm{co}(k,m) \wedge k \leq m\}=$$
$$ \{(k,n);\ \textrm{co}(k,n) \wedge k \leq n\}.$$

\begin{itemize}
\item $\mathbb I_1 = \{(1,1)\}, \mathbb I_2 = \{(1,2)\}, \mathbb I_3 = \{(1,3), (2,3)\}, \mathbb I_4 = \{(1,4), (3,4)\}$, 

$\mathbb I_5 = \{(1,5), (2,5), (3,5), (4,5)\}, \mathbb I_6 = \{(1,6), (5,6)\},  \dots$

\item The {characteristic sequence} 
$$\chi(\mathbb I) = (1,1,2,2,4,2,6,4, \dots).$$ 
The $n$-th term of the characteristic sequence is equal to the number of coprime numbers less than $n$ that is given by \emph{Euler's function $\varphi(n)$.}
$$\chi_n(\mathbb I) = \varphi(n) \leq n-1.$$
\item The {size sequence} 
$$\sigma(\mathbb I) = (1,2,4,6,10,12,18,22,28, \dots).$$ 
The $n$-th term of the size sequence
$$ \sigma_n(\mathbb I) = \chi_1(\mathbb I) + \dots + \chi_n(\mathbb I) = \sum _{i=1}^{n}\varphi(i) = \Phi(n).$$
This sum is the \emph{totient summatory function} usually denoted by $\Phi(n)$. 
\end{itemize}

\begin{theorem}\label{unit} The unit interval \ $\mathbb I = (0,1]_\mathbb Q$ of rational numbers is canonically countable. Its size $\sigma(\mathbb I)$ is given by the totient summatory function $\Phi$ and is algorithmically enumerable. It is equal in order to $\alpha^2$.
$$\sigma(\mathbb I) = (\Phi(n))_{n\in \mathbb N} \approx_\mathcal F \alpha^2, $$
More precisely 
$$\frac{3}{10} \cdot \alpha^2 < \sigma(\mathbb I) < \frac{\alpha^2 - \alpha}{2}$$
\end{theorem}
\begin{proof} According to Theorem 330 in 
(Hardy \& Wright 1968, p. 268), it is valid for the totient summatory function $\Phi$: $$\Phi(n) =\frac{3 n^2}{\pi^2} + O((\log n)^{2 \over 3}(\log \log n)^{4 \over 3}) .$$ Thus, and also because $\varphi(i) \leq i-1$ for $i \in \mathbb N$
$$\frac{3}{10} \cdot n^2< \frac{3 n^2}{\pi^2} < \Phi(n) = \sigma_n(\mathbb I) = \sum _{i=1}^{n}\varphi(i) \leq \sum _{i=1}^{n}i -1 = \frac{n \cdot (n - 1)}{2}.$$ 
\end{proof}


\begin{definition} We denote the set size of the half-open unit interval $\mathbb I = (0,1]_\mathbb Q$ of rational numbers by the letter $\varphi$. 
$$ \varphi = \sigma(\mathbb I) = (1,2,4,6,10,12,18,22,28, \dots). $$
\end{definition}

The result that $\sigma(\mathbb I) \approx \alpha^2 $ may be a little surprising at first sight, but the reasoning behind it is quite intuitive. The character of rational numbers in the interval is different from that of natural numbers. The infinite set of the former is literally denser; between any two rational numbers there is another, not just one, but infinitely many rational numbers. Unlike the natural numbers, which are ordered in a sequence where each number has an immediate successor. 

\subsubsection{The Problem of Rational Numbers}

Although it might seem appropriate to represent all positive rational numbers $\mathbb Q^+$ as coprime pairs of natural numbers we show that this is not the right way. Let 
$$\mathbb Q^+ = \{\frac{k}{m} \sim (k,m) \in \mathbb N \times \mathbb N; \textrm{co}(k,m)\}. $$

Then the $n$-th component of this arrangement would be
$$\mathbb Q^+_n = \mathbb Q^+ \ \cap \ (\mathbb N \times \mathbb N)_n = \{(k,m) \in \mathbb N \times \mathbb N; \textrm{co}(k,m) \wedge n = \max\{k,m\}\}$$
$$\mathbb Q^+_1 = \{(1,1)\}, \mathbb Q^+_2 = \{(1,2), (2,1)\}, \mathbb Q^+_3 =\{(1,3),(2,3),(3,2),(3,1)\}, \dots$$ 
This arrangement is an extension of the canonical arrangement of both natural numbers and the unit interval of rational numbers. So far, it is correct. The problem lies elsewhere. 
Positive rational numbers greater than or equal to $1$ are contained in the interval 
$$\mathbb J =  [1, \infty)_\mathbb Q = \{\frac{k}{m} \sim (k,m) \in \mathbb N \times \mathbb N; \textrm{co}(k,m) \wedge k \geq m\} $$
The $n$-th component of the canonical arrangement would be 
$$\mathbb J_n = \mathbb J \ \cap \ (\mathbb N \times \mathbb N)_n = \{(n,m);\ \textrm{co}(n,m) \wedge n \geq m\}.$$
 $$\mathbb J_1 = \{(1,1)\}, \mathbb J_2 =  \{(2,1),\}, \mathbb J_3 = \{(3,1),(3,2)\}, \mathbb J_4 = \{(4,3),(4,1)\}, \dots $$ 
$$\chi(\mathbb J) = (1,1,2,2,4, \dots) = \chi(\mathbb I)$$
The size sequences of $\mathbb J$ would be the same as that of $\mathbb I$.  
$$\sigma(\mathbb J) = (1,2,4,6,10, \dots) = \sigma(\mathbb I)$$
The size of all positive rational numbers $\mathbb Q^+$ would be twice as great as the size of $\mathbb I$ minus $1$.
$$\sigma(\mathbb Q^+)  = \sigma(\mathbb I \cup \mathbb J) = \sigma(\mathbb I) + \sigma(\mathbb J) - \sigma(\mathbb I \cap \mathbb J) = 2 \cdot \sigma(\mathbb I) - 1.$$ 
 That is not satisfactory. Among others, this representation is not homogeneous, it does not preserve the homogeneity of size. contradicts Bolzano's requirement of homogeneity. 
A much better option seems to be to express positive rational numbers as mixed fractions.

\subsubsection{Positive Rational Numbers}
 Every positive rational number can be expressed as a mixed fraction $p + {k \over m}$ where $p \in \mathbb N_0$ and $k,m \in \mathbb N$ are coprime. It is represented by the triple $(p,k,m)$.\footnote{Of course, the triple $(p,1,1)$ represents $p + 1$. } Thus, the positive rational numbers are represented as a Cartesian product of $\mathbb N_0$ and $\mathbb I$.

\begin{definition} 
We introduce a predicate $\textrm{fr}(p,k,m)$ saying that the triple $(p,k,m) \in \mathbb N_0 \times \mathbb N^2$ represents a positive mixed fraction $p + {k \over m}$
$$\textrm{fr}(p,k,m) \Leftrightarrow ((p,k,m) \in \mathbb N_0 \times \mathbb N^2 \wedge  \textrm{co}(k,m)\wedge k \leq m)$$
The positive rational numbers are represented by the set 
$$\mathbb Q^+  = \mathbb N_0 \times \mathbb I= \{(p,k,m); \textrm{fr}(p,k,m)\}.$$
\end{definition}

Both $ \mathbb N_0$ and $\mathbb I$ are canonically arranged as well as their Cartesian product. 
So, the canonical arrangement of rational numbers is determined. Its $n$-th component   
$$\mathbb Q_n^+ = \mathbb Q^+ \cap \ ({\mathbb N_0}^3)_n$$ 
\begin{itemize}
\item $\mathbb Q_n^+ =  \{(p,k,m); \textrm{fr}(p,k,m) \wedge n = \max\{p,k,m\}\} =$

$\quad  \quad \ \ \  \{(p,k,m); \textrm{fr}(p,k,m) \wedge n = \max\{p,m\}\} =$


$\quad \quad \ \ \  \{(p,k,m); \textrm{fr}(p,k,m) \wedge ((p = n \wedge m \leq n ) \vee ( m = n \wedge p < n ))\}. $


\item $\mathbb Q^+_1 = \{(0,1,1),(1,1,1),\},$
 
$\mathbb Q^+_2 =\{(0,1,2),(1,1,2),(2,1,1),(2,1,2)\}$, 

$\mathbb Q^+_3=\{(0,1,3), (0,2,3), (1,1,3), (1,2,3), (2,1,3), (2,2,3), (3,1,1),$ 

$(3,1,3), (3,1,2), (3,2,3)\}, \dots$
\item The characteristic sequence $\chi(\mathbb Q^+) = (2,4, 10,14,30, \dots)$.
\item The size sequence
$\sigma(\mathbb Q^+) =  (2,6,16,30,60 \dots ).$ 
\end{itemize}
\begin{theorem} Positive rational numbers are canonically countable. Their size is 
$$\sigma(\mathbb Q^+) = (\alpha +1) \cdot \varphi \approx_\mathcal F \alpha^3.$$ 
More precisely 
$${3 \over 10}(\alpha^3 + \alpha^2) <_\mathcal F \sigma(\mathbb Q^+) <_\mathcal F \frac{ \alpha^3 - \alpha}{2}.$$ 
\end{theorem}
\begin{proof} 
This arrangement of rational numbers is the extension of the canonical arrangements both of natural numbers and of the unit interval of rational numbers.
 $$\sigma(\mathbb Q^+) = \sigma(\mathbb N_0) \cdot \sigma(\mathbb I) = ( \alpha +1) \cdot \varphi =
(2,3,4, 5\dots) \cdot (1,2,4,6 \dots) = (2,6,16,30\dots )$$
According to Theorem \ref{unit}
$$ (\alpha +1) \cdot {3 \over 10}\alpha^2 <_\mathcal F \sigma(\mathbb Q^+)= \sigma(\mathbb N_0) \cdot \sigma(\mathbb I)  <_\mathcal F (\alpha+1) \cdot \frac{\alpha^2 - \alpha}{2}= \frac{\alpha^3 - \alpha}{2}.$$
\end{proof}

\begin{theorem}\label{rat}
Rational numbers are canonically countable. Their size is
$$\sigma(\mathbb Q) = 2  \varphi \cdot (\alpha + 1) - 1 \approx_\mathcal F \alpha^3.$$
More precisely
$${3 \over 5}(\alpha^3 + \alpha^2)  <_\mathcal F \sigma(\mathbb Q) <_\mathcal F \alpha^3 - \alpha.$$
\end{theorem}
\begin{proof}
Rational numbers $\mathbb Q$ are the union of positive and negative rational numbers and zero. 
$$\mathbb Q = \mathbb Q^+ \cup \mathbb Q^- \cup \{0\}.$$ We assume that negative rational numbers have the same size as positive ones. So
$$2 \cdot {3 \over 10}(\alpha^3 + \alpha^2) + 1 <_\mathcal F \sigma(\mathbb Q) = 2 \cdot \sigma(\mathbb Q^+) + 1 <_\mathcal F 2 \cdot \frac{ \alpha^3 - \alpha}{2} + 1.$$ 
\end{proof}

We proved that $\sigma(\mathbb Q) \approx \alpha^3$, more precisely ${3 \over 5}(\alpha^3 + \alpha^2)  <_\mathcal F \sigma(\mathbb Q) <_\mathcal F \alpha^3 - \alpha.$ 
The intuitive explanation is that the rational numbers are infinite not only \enquote{in length} but also \enquote{in depth}. The nature of the two infinities is different. The set of rational numbers between any two successive integers is infinitely denser than the set of integers; its size is infinitely greater than the size of integers. Rational numbers are the union of infinitely many such intervals.   

Of course, the result depends on the representation of rational numbers. But this representation, as   %
the Cartesian product $\mathbb N \times \mathbb I$, is quite intuitive, 
it depends on a determining ground. It extends the representations both of $\mathbb N$ and of $\mathbb I$. 
Moreover, it is \emph{homogeneous}, two intervals of rational numbers of the same length have the same size, as the following theorem demonstrates. 
\begin{theorem}\label{interval} All intervals of rational numbers of the same length have the same size.
\end{theorem}
\begin{proof}  We prove it in several steps. Let $j,k,m,n,p,q \in \mathbb N$.
\begin{enumerate}[(1)]
\item \emph{Every unit interval starting at $j \in \mathbb N$ has the same size as $\mathbb I$}: 
 
Let us denote $\mathbb J = (j,j+1]_\mathbb Q$. Then its canonical arrangement 
$$\mathbb J_n = \mathbb Q^+_n \cap (j,j+1]_\mathbb Q $$
From now on, we shall assume for all pairs $(k,m)$ that $\textrm{co}(k, m) \wedge k \leq m$.
$$\mathbb J_n =  \{(p,k,m); \textrm{fr}(p,k,m) \wedge n = \max\{p,m\}\}\cap (j,j+1]_\mathbb Q  = $$
$$\{(j,k,m); \textrm{fr}(j,k,m) \wedge n = \max\{j,m\}\} = $$
$$\{(j, k, m);  \textrm{fr}(j,k,m) \wedge ( n=j  \wedge m \leq n ) \vee (n =m \wedge j \leq n )\}$$
There are three options:
\begin{itemize}
\item If $n < j$ then $\mathbb J_n = \emptyset$, $\chi_n(\mathbb J) = 0$.
\item If $n= j$ then $\mathbb J_n = \{(j,k,m);  \textrm{fr}(j,k,m) \wedge m \leq n\}$, $|\mathbb J_n| = \chi_n(\mathbb J) = \sigma_n(\mathbb I)$.
\item If $n > j$ then $\mathbb J_n = \{(j,k,m);  \textrm{fr}(j,k,m) \wedge m = n\}$, $|\mathbb J_n| = \chi_n(\mathbb J) = \chi_n( \mathbb I)$. 
\end{itemize}
Consequently, the size sequence 
$$\sigma(\mathbb J)_\mathbb Q =_\mathcal F \sigma(\mathbb I).$$

\item \emph{Every unit interval $\mathbb J$ starting at any rational number has the same size as $\mathbb I$.} 

Let $\mathbb J = (j+\frac{p}{q},(j+1) +\frac{p}{q}]_\mathbb Q$. 
We will not write all details. We again get three similar options. The only difference is when $n= j$. However, $\chi_n(\mathbb J) + \chi_{n+1}(\mathbb J) = \sigma_{n+1}(\mathbb I)$. The result is the same  
$$\sigma(j+\frac{p}{q}, (j+1) +\frac{p}{q}]_\mathbb Q =_\mathcal F \sigma(\mathbb I).$$

\item \emph{Every two intervals of the same length has the same size.}

The proof proceeds in a similar way, but it is technically more complicated. First, we prove that the interval $(0, {p \over q}]_\mathbb Q$ has the same size as $(r, r + {p \over q}]_\mathbb Q$ for $r \in \mathbb N$. Then we shift the start of the interval to any rational number $r \in \mathbb Q$.  
\end{enumerate}
\end{proof}


\section{Numerosity Theory}\label{num}

The present theory of canonically countable sets and their sizes is similar to \emph{Numerosity Theory} (NT), which has been developed since 1995, first formalized in (Benci \& Di Nasso 2003), and published in a book (Benci \& Di Nasso 2019).\footnote{It might seem that we could call our theory a \enquote{cheap version of numerosity theory}, similar to how Terence Tao describes a \enquote{cheap version of non-standard analysis} (Tao 2012a). However, that would discredit it. 
Tao's theory is an interesting attempt to reconstruct non-standard analysis using only the Fr\' echet filter instead of an ultrafilter. It is constructive, but of course less powerful. The structure $\mathbb R^\mathbb N/\mathcal F$ contains unnecessary elements and non-zero divisors of zero. The transfer principle is not valid. We also deal with non-decreasing sequences of natural numbers factorized by the Fr\' echet filter instead of an ultrafilter. For our purposes, however, this structure is sufficient. It does not contain unnecessary elements and the transfer principle is not needed.} NT is a part of Alpha-Theory, which is a new elementary axiomatics which can serve as a basis for infinitesimal calculus. 
We appreciate this theory and are grateful for its solid mathematical foundation. However, our aim is a bit different. We are primarily concerned with the sizes of sets and do not intend to provide a basis for non-standard analysis as NT also does. 

\subsection{Basic Definitions}

The key notion of NT is the \emph{labelled} set which is a pair $(A, l)$ where $A$ is a set and $l: A \longrightarrow \mathbb N_0$ is a \emph{finite-to-one} function which means that all pre-images $l^{-1}(n)$ for all $n \in \mathbb N$ are finite. The function $l$ is called a \emph{labelling function} and the natural number $l(x)$ is a \emph{label} of $x$.

A canonically countable set $A = \bigcup\{A_n; n \in \mathbb N\}$ in our theory is naturally labelled. The function $l$ assigns to every element of $A$ the index of its component. 
$$x \in A_n \Leftrightarrow l(x) = n.$$

Labelled sets are \emph{ex definitione} countable, but they need not be canonically countable. Although NT speaks about \emph{canonical labelling} of natural numbers and   integers which corresponds to ours, the canonical arrangement is not considered substantial. We define sizes of sets, their union and the Cartesian product on the basis of the canonical  arrangement, NT simply defines these notions without providing any justification. However, both theories lead to similar results (Benci \& Di Nasso 2019, p. 278 -- 280). 

The \emph{counting function} $\varphi_A: \mathbb N_0 \longrightarrow \mathbb N_0 $ is defined 
$$\varphi_A(n) = |\{a \in A; l(a) \leq n\}|$$

The \emph{approximating sequence}\footnote{(Benci \& Di Nasso 2003, p. 52).}  $(\varphi_A(n))_n$ is  equal to our size sequence $\sigma(A)$. 

The \emph{$\alpha$-numerosity} of $A$ is defined 

$$\mathfrak n_\alpha(A) = \lim_{n\uparrow \alpha}{\varphi_A (n)}$$

In the basic model of Alpha-Theory, $\lim_{n\uparrow \alpha}{\varphi_A (n)} = [\varphi]_\mathcal U$ is an equivalence-class of the ultraproduct $\mathbb R^\ast = \mathbb R^\mathbb N/ \mathcal U$ where $\mathcal U$ is a non-principal ultrafilter on $\mathbb N$. (Benci \& Di Nasso, 2019, p. 221). 

\subsection{Ultrafilters versus the Fr\' echet Filter}\label{ultr}

As mentioned in Section \ref{partial order}, the greatest objection against our theory might be that the set sizes are not linearly ordered. The linear ordering of numerosities is enabled just by an ultrafilter. If we used a non-principal ultrafilter $\mathcal U$ instead of mere Fr\' echet filter $\mathcal F$ in Definition \ref{Bdef} we would obtain the same model as NT. 

\begin{definition} 
 Let $(a_n), (b_n)$ be two sequences of natural numbers. 
$$(a_n) =_\mathcal U (b_n) \text{ if and only if } \{n; a_n = b_n \} \in \mathcal U.$$
$$(a_n) <_\mathcal U (b_n) \text{ if and only if }  \{n; a_n <  b_n \} \in \mathcal U.$$
\end{definition}
The equivalence-classes $[(a_n)] \in \mathbb N^\ast = \mathbb N^\mathbb N/ \mathcal U$ of the ultraproduct are linearly ordered. The utrafilter is a powerfull instruments. Thanks to its wonderful properties, the infinitesimal calculus has been consistently formalized. However, it is a non-constructive object, whose existence is guaranteed by the \emph{Axiom of Choice}. One never knows all its properties. Consequently, some of the results obtained with an ultrafilter are arbitrary.\footnote{The objection against the arbitrariness arising from the use of an ultrafilter in NT does not apply to the use of an ultrafilter in non-standard analysis and other hyperreal models where it is a very useful, perhaps even necessary instrument. 
The alleged arbitrariness is rebutted in (Botazzi \& Katz 2021) by showing, among other things, that their critics ignore the property of being internal. If one considers only hyperfinite measures, no underdetermination occurs; the internal results of functions and measures do not depend on the choice of an ultrafilter.

The situation for numerosities is different. The internal results such as whether there are fewer or as many even as odd numbers or whether $\alpha$ is odd or even depend directly on the choice of an ultrafilter.} To the credit of NT, it should be said that it is aware of this problem and highlights it. 
\begin{quote} We tell in advance that every model of numerosity is grounded on an ultrafilter; in consequence, many features of the theory will depend on the choice of such an underlying ultrafilter. (Benci \& Di Nasso 2019, p. 292). 
\end{quote}

It is partially possible to get around this problem. In (Benci \& Di Nasso 2019, p. 288 - 289), the authors take a selective ultrafilter containing sets of all $k$-multiples $\mathbb M_k = \{k, 2k, 3k, \dots \}$ and also sets of all $k$-th powers $\mathbb S_k = \{1^k, 2^k, 3^k, \dots \}$ of all $k \in \mathbb N$. This is conceivable because these sets form a centered system that can be extended to an ultrafilter. 

Then, the infinite number $\alpha$ is divisible by any $k \in \mathbb N$ since it is the product of $k$ and the numerosity of the set of $k$-multiples. Moreover, $\alpha$ is a $k$-th power of the numerosity of the set of $k$-th powers of any $k \in \mathbb N$ . 
$$\alpha =_\mathcal U k \cdot \mathfrak n(\mathbb M_k) .$$
$$\alpha =_\mathcal U  \mathfrak n(\mathbb S_k)^k.$$

One can raise several objections against these equations. First, they are valid just for the chosen ultrafilter. Second, there are still uncountably many ultrafilters with the required properties, whereas the other properties remain arbitrary. Third, the only reason to take this particular ultrafilter is simplicity. The same reason could be given for the validity of the continuum hypothesis, but this is not considered sufficient.

What we can say for granted is that there is a consistent model of NT in which these equations hold. We can equally say there is a consistent model of NT in which none of them holds, for instance if we take a selective ultrafilter containing sets $\{k+1, 2k+1, 3k+1, \dots \}$ and also sets  $\{1^k+1, 2^k+1, 3^k+1, \dots \}$ for all $k \in \mathbb N$. It is also conceivable for the same reason. 

We \emph{proved} a similar but a little weaker Theorem \ref{prop} saying that the difference between $\alpha$ and the product of $k$ and the size of the set of $k$-multiples is less than or equal to $k$. This assertion would be valid even if we eventually used any ultrafilter. 
$$\alpha - k \leq_\mathcal F k \cdot \sigma((\mathbb M_k) \leq_\mathcal F \alpha$$

NT has one more reason for the need of a linear ordering of numerosities and so of an ultrafilter. Numerosities are used as a basis for the Alpha-Calculus that is a special kind of non-standard analysis. 
\begin{quote} \dots because it is the very idea of numerosity that leads to Alpha-Calculus. (Benci \& Di Nasso 2019, p. 300) \end{quote}

We do not intend to use size sequences for this purpose or anything like that. All we want is to determine the sizes of countable sets so that PW holds and cardinalities of finite sets are preserved. We believe that the Fr\'{e}chet filter is sufficient for this purpose.

\subsection{The Case of Rational Numbers}

A canonical arrangement of rational numbers is a more delicate question. We have explored this point in Section \ref{rational}
and considered the unit interval $\mathbb I = (0, 1]$ as a subset of  $\mathbb N \times \mathbb N$ of coprime numbers such that denominator is less than numerator. 
Rational numbers $\mathbb Q$ are represented as a subset of $\mathbb Z \times \mathbb I$, their size is 
 $$\sigma(\mathbb Q) =  2  \varphi \cdot (\alpha + 1) - 1 \approx_\mathcal F \alpha^3.$$
 
In contrast, NT claims:
\begin{quote} In particular, it seems there is no definitive way to decide whether $\mathfrak n_\alpha((0,1]_\mathbb Q) \geq \alpha$ or $\mathfrak n_\alpha((0,1]_\mathbb Q) \leq \alpha$. So, in absence of any reason to choose one of the two possibilities, we go for the simplest option $\mathfrak n_\alpha((0,1]_\mathbb Q) = \alpha$. (Benci \& Di Nasso 2019, p. 291) 
\end{quote}

Consequently, the numerosity of positive rational numbers $\mathfrak n(\mathbb Q^+) = \alpha^2$ and the numerosity of all rational numbers $$\mathfrak n(\mathbb Q^+) = 2 \alpha^2 +1.$$ 
The difference in the results is due to the distinct representations of the rational numbers. From the beginning, we have justified every step we have taken so as not to drift into the arbitrariness.  
To support the soundness of our concept, we point out Theorem \ref{interval} that proves that all intervals of the same length have the same size.

\section{Conclusion}\label{concl}

In his paper \emph{Set Size and the Part-Whole Principle}, Parker argues that there can be no good theories of set size satisfying PW, which he calls Euclidean theories. 

\begin{quote}But our main question here is whether it is possible to have a really good Euclidean theory of set size, and my answer is, no, not really - not if that means it must be strong, general, well-motivated, and informative. I will argue that any Euclidean theory strong
and general enough to determine the sizes of certain simple, countably infinite sets must
incorporate thoroughly arbitrary choices.\footnote{Parker explains what he means by these terms in a footnote 5 on the same page. \enquote{By a general theory I mean one that applies to a broad domain or many domains, including,
in particular, subsets of the whole numbers and countable point sets. By strong I mean logically
strong—a theory that leaves little undecided, and in particular, one that determines the sizes of
the simple sets discussed below. By well motivated I mean that the details of the
theory—all of the particular sizes and size relations it assigns—are so assigned for some reason;
they are not chosen arbitrarily. And finally, informative here means that the consequences of the
theory indicate something of interest that holds independently of the theory itself.}} (Parker 2013, p. 590). \end{quote}

The present theory is motivated by Bolzano's principles: the canonical arrangement of sets is based on his notion of a \enquote{determining ground}, the size sequences and their structure are derived from Bolzano's infinite series, they generalize the notion of finite cardinalities.\footnote{Although there is a notion of the asymptotic density, as used in number theory, that serves as a mathematical tool for discriminating sizes of infinite sets of natural numbers, it is not sufficient. It does not generalize the notion of cardinality of finite sets; all finite sets have the same density, namely 0. Some infinite sets have no density and some infinite sets have density 0. PW is respected only partially. (Mancosu 2009, p. 627)} PW is valid up to one element. The sizes of integers, rational numbers, their subsets and their products are uniquely determined and algorithmically enumerable. Moreover, we proved that the size of rational numbers is homogeneous. These arguments are perhaps sufficient to show that this theory is \enquote{good} in Parker's sense and is not \enquote{arbitrary and misleading}.  

Another thing is that when Parker investigated the possible existence of theories preserving PW he assumed that an assignment determining the size of sets is a function from a class of sets $D$ to a \emph{linearly} ordered mathematical structure.\footnote{However, Parker himself 
admits that sizes of some sets may be incomparable.\enquote{An assignment is Euclidean on $D$ if PW applies to all proper subset/superset pairs in $D$. In that case, all such pairs have size relation, but this does not imply totality; we might for example, have disjoint sets in $D$ that are not comparable at all.} (Parker 2013, p. 593).}  That is not our theory. The codomain of $\sigma$  is the set of non-decreasing sequences of natural numbers modulo the Fr\' echet filter which is just partially and not linearly ordered. 

The Fr\' echet filter $\mathcal F$ means a simple, constructive and uniquely defined criterium for order and equality of two size sequences - order or equality is valid starting from a sufficiently great index of terms. It is contained in any ultrafilter defined on $\mathbb N$,  $\mathcal F \subseteq \mathcal U$. The results are ultimate and would be valid in \emph{any} eventual extension of the theory via the ultrafilter. And so, the use of the Fr\'{e}chet filter determines the limits of our certain knowledge. Everything beyond them is hypothetical. 

On the other side, the ordering of set-sizes is indeed really just partial which means that some sizes are not comparable and some properties of infinite sets cannot be determined, such as whether $\alpha$ is even or odd.\footnote{Benci and Di Nasso (2019, p. 275) also point out that the \emph{trichotomic property} of equisize relation, which implies linearity of set sizes, is disputable. Zermelo needed the full strength of the \emph{axiom of choice} to prove it for Cantorian cardinalities.} \emph{Quid pro quo.}



%

Cantorian cardinalities and set sizes are not mutually exclusive concepts. They are just two different points of view. According to Hume's Principle, the natural numbers and the squares have the \emph{same cardinality}. However, they have \emph{different sizes} for which PW applies. Galileo's paradox is just a puzzle with two possible solutions depending on the \enquote{determining ground} of the squares.  \medskip

\noindent \footnotesize \textbf{Acknowledgement.} I would like to thank to Petr K\accent23 urka, Piotr B\l{}aszczyk and Ansten Klev for helpful comments on an earlier version of the manuscript. 
I am very grateful to my anonymous referee for the thorough reading of the manuscript and valuable comments.   


\begin{thebibliography}{}
\bibitem{} Bair, J., B\l{}aszczyk, P., Ely, R., Henry, V., Kanovei, V., Katz, K., Katz, M., Kutateladze, S., McGaffey, T., Schaps, D., Sherry, D., Shnider, S. (2013). Is Mathematical History Written by the Victors? \emph{Notices of the American Mathematical Society} 60(7), 886-904.  
\bibitem{} Bartlett, J., (2021). The Products of Hyperreal Series and the Limitations of Cauchy Products. \emph{Communications of the Blyth Institute} 3(2), 34 - 36.
\bibitem{} Bartlett, J., Logan G. \& Nemati, D., (2020). Hyperreal Numbers for Infinite Divergent Series. \emph{Communications of the Blyth Institute} 2(1), 7 - 15.
\bibitem{} Bellomo A., \& Massas, G., (2021). Bolzano's Mathematical Infinite. \emph{The Review of Symbolic Logic} 17(1), 1 - 55.
\bibitem{} Benci V. \& Di Nasso, M., (2003). Numerosities of Labelled Sets: a New Way of Counting. \emph{Advances in Mathematics} 173, 50 - 67.
\bibitem{} Benci, V. \& Di Nasso, M., (2019). \emph{How to Measure the Infinite, Mathematics with Infinite and Infinitesimal Numbers}, World Scientific.
\bibitem{} B\l{}asczyk, P., (2021). Galileo's Paradox and Numerosities, \emph{Philosophical Problems in Science (Zagadnienia Filozoficzne W Nauce)} 70, 73 - 107.
\bibitem {} Bolzano, B., (1851/2004). \emph{Paradoxien des Unendlichen}, CH Reclam, Leipzig. English translation \emph{Paradoxes of the Infinite} in (Russ 2004).
\bibitem{} Bottazzi, E. \& Katz, M., (2021). Infinite Lotteries, Spinners, Applicability of Hyperreals, \emph{Philosophia Mathematica} 29(1), 88–109.
\bibitem {} Cantor, G., (1883/1976). \emph{Grundlagen einer allgemeinen Mannigfaltigkeitslehre}. English translation Foundations of a General Theory of Manifolds by Georg Cantor, \emph{The Campaigner, Journal of the National Caucus of Labor Committees,} 9(1-2), 69 - 97. 
\bibitem {} Ferreiros, J., (1999). \emph{Labyrinth of Thought. A History of Set theory and its Role in Modern Mathematics}. Birkh\H auser Verlag AG, Basel - Boston - Berlin. 
\bibitem{} Galileo, G., (1638/1914). \emph{Dialogues Concerning Two New Sciences}. English translation Henry Crew \& Alfonso da Savie, New York, Macmillan, http://files.libertyfund.org/files/753/Galileo 0416 EBk v6.0.pdf.
\bibitem{} Guill\' en, E.F., (2021). Bolzano's Theory of messbare Zahlen: Insights and Uncertainties Regarding the Number Continuum. In: Sriraman, B. (eds) \emph{Handbook of the History and Philosophy of Mathematical Practice}, 1 - 38,  Springer, Cham. 

\bibitem{} Hardy, G.H. \& Wright, E.M., (1968). \emph{The Introduction to the Theory of Numbers}. Oxford University Press. Ely House, London.
\bibitem{} Heck, R.K., (2019) The Basic Laws of Cardinal Number.
In Philip A. Ebert \& Marcus Rossberg (eds.), \emph{Essays on Frege's Basic Laws of Arithmetic}. Oxford University Press, 1 - 30. 
\bibitem{} Jech, T., (2006). \emph{Set Theory}. Springer Verlag. Berlin - Heidelberg - New York. 
\bibitem{} Kanamori, A., (2012). In Praise of Replacement. \emph{The Bulletin of Symbolic Logic} 18(1), 46 - 90.
\bibitem{} Knuth, D., (1976). Big Omicron and big Omega and big Theta. \emph{ACM SIGACT News} 8(2), 18 - 24. 
\bibitem {} Mancosu, P., (2009). Measuring the Size of Infinite Collections of Natural Numbers: Was Cantor's Set Theory Inevitable? \emph{The Review of Symbolic Logic} 2(4), 612 - 646. 
\bibitem{} Mancosu, P., (2015). In Good Company? On Hume's Principle and the Assignement of Numbers to Infinite Concepts. \emph{The Review of Symbolic Logic} 8(2),  370 - 410. 
\bibitem{} Parker, M.W., (2009). Philosophical method and Galileo's paradox of infinity. In \emph{New Perspectives on Mathematical Practices}. Hackensack, NJ, World Scientific, 76 - 113. 
\bibitem{} Parker, M.W., (2013). Set-size and the Part-Whole Principle. \emph{The Review of Symbolic Logic} 6(4), 589 - 612.
\bibitem{} Rosser, J.B. \& Schoenfeld, L., (1962). Approximate Formulas for some Functions of Prime Numbers. \emph{Illinois J. Math.} 6(1), 64 - 94. 
\bibitem{} Rusnock, P. \& \v Sebest\' ik, J., (2019). \emph{Bolzano: His Life and Work.} Oxford University Press. 
\bibitem {} Russ S., (2004). \emph{The Mathematical Works of Bernard Bolzano}. Oxford University Press.
\bibitem {} Russ, S. \& Trlifajov\'{a}, K., (2016). Bolzano's Measurable Numbers: are they Real? In \emph{Research in History and Philosophy of Mathematics}. Basel, Birkh\H{a}user, 39 - 56. 
\bibitem {} \v{S}ebest\'{i}k, J., (1992). \emph{Logique et math\'{e}matique chez Bernard Bolzano}. Vrin, Paris.
\bibitem{} Simon, P., (1998). Bolzano on Collections. In \emph{Bolzano and Analytic Philosophy}, ed. by W. K\H{u}nne, M. Siebel and M. Textor, Brill, Rodopi.
\bibitem {} Tao, T., (2012). \emph{A Cheap Version of Non-standard Analysis}, 

https://terrytao.wordpress.com/2012/04/02/a-cheap-version-of-nonstandard-analysis/. 
\bibitem{} Trlifajov\'{a}, K., (2018). Bolzano's Infinite Quantities, \emph{Foundations of Science} 23(4), 681 - 704.

\end{thebibliography}
\end{document}